\documentclass[twoside,regno,12pt]{article}
\usepackage{amssymb}
\newcommand{\rad}{\mathop{\rm rad}}
\begin{document}

\begin{center}

{\bf  THE VARIATION   OF THE RANDI\' C  INDEX \\
WITH REGARD TO MINIMUM AND MAXIMUM DEGREE  }\\
\vskip 0.4cm
 Milica Milivojevi\'c,Ljiljana Pavlovi\'c\footnote{corresponding author ( pavlovic@kg.ac.rs)}, \\
 \vskip 0.2cm
{\small\it Faculty of Science and Mathematics, University of Kragujevac,}\\
{\small\it Radoja Domanovi\'ca 12, Kragujevac, Serbia.}\\
\end{center}
\vskip\baselineskip
\par\noindent {\bf Abstract.}
The variation of the Randi\'c index $ R'(G) $ of a graph $G$ is defined
by\ $R(G) = \sum_{uv \in E(G)}\frac 1{\max \{d(u) d(v)\}}$, where
$d(u)$ is the degree of vertex $u$ and the summation extends
over all edges $uv$ of $G$. Let $G(k,n)$ be the set of connected simple  $n$-vertex graphs
with minimum  vertex degree $k$. In this paper we  found
in $G(k,n)$ graphs for which the variation of the Randi\'c index
 attains its minimum value.
 When $k \leq \frac n2$  the extremal graphs are complete split graphs $K_{k,n-k}^*$, which   have
 only vertices of two degrees, i.e. degree $k$ and degree $n-1$, and the number of vertices of degree $k$
 is $n-k$, while  the number of vertices of degree $n-1$ is $k$.  For $k \geq \frac n2$  the extremal graphs
have also  vertices of two degrees $k$ and $n-1$, and the number of vertices of degree $k$ is $\frac n2$.
Further, we generalized results for graphs with given maximum degree.

\par {\bf Keywords:}  Simple graphs with given minimum degree,
 Variation of the Randi\'c index, Combinatorial optimization, Quadratic programming.

\par {\bf AMS Subject classifications:} 05C35, 90C27.

 \vskip\baselineskip
\noindent{\bf  1. INTRODUCTION}
\vskip\baselineskip

In 1975 Randi\'c
proposed a topological index, suitable for measuring the extent of
branching of the carbon-atom skeleton of saturated hydrocarbons.
The Randi\'c index $R(G)$ of a graph $G$, defined in
[13], is given by
$$
R(G) = \sum_{uv \in E(G)}\frac 1{\sqrt {d(u) d(v)}},
$$
where the summation extends over all edges  of $G$ and $d(u)$ is the degree of the vertex $u$ in $G$. Randi\'c  himself
demonstrated [13] that this index is well correlated with a variety
of physico-chemical properties of alkanes.
 The Randi\'c index has  become one of the most popular molecular descriptors.
 To this index several books are devoted ([8-10]).
 Later, in 1998 Bollob\'as  and Erd\H{o}s [3] introduced general Randi\'c index $R_{\alpha}$, where $\alpha$ ia a real number, as
$$
R_{\alpha}(G) = \sum_{uv \in E(G)}( d(u) d(v))^{\alpha}.
$$

In order to attack some conjectures concerning the Randi\'c index,
Dvo\v{r}\'ak et al.
% [European J. Combin. 32 (2011), 434-442]
introduced in [6 ] a variation of this index, denoted by $R'$.
The variation of the Randi\'c index of  a graph $G$ is given by
$$
R'(G) = \sum_{uv \in E(G)}\frac 1{\max \{d(u) d(v)\}}.
$$
Although no application of the $R'$ index in chemistry is known so far, nevertheless this index turns out to be very useful, especially from a mathematical point of view, since it is much easier to follow  during graph modifications than the Randi\'c index.  Using the $R'$ index, Cygan et al. [4]
% proved that for any connected graph $G$ of maximum degree at most four which is not a path with even number of vertices, $R(G) \geq \rad(G)$. As a consequence, they
resolved the conjecture $R(G) \geq \rad(G) -1$ given by Fajtlowicz 1988 in [7] for the case when $G$ is a chemical graph.
% They actually showed that for all connected chemical graphs $G$ the inequality
% $R'(G) \geq \rad(G) -\frac 12$ holds.
In [1] Andova et al.
 determined graphs with minimal and maximal value for  the $R'$  index , as well as graphs with minimal and maximal value of the $R'$ index among trees and unicyclic graphs.
They also showed that if $G$ is a triangle free graph on $n$ vertices with minimum degree  $\delta (G)$, then $R'(G) \geq \delta$.

\par Now we define terms and symbols used in the paper.
Let  $G(k,n)$ be the  set of connected simple $n$-vertex graphs
with minimum  vertex degree $k$. If  $u$ is a vertex of $G$, then
 $d(u)$ denotes the degree of the vertex $u$, that is, the
number of edges of which $u$  is an endpoint.
Let  $V(G)$, $E(G)$,  $\delta (G)$  and $\Delta(G)$  denote the vertex set, edge set, minimum degree, and maximum degree of $G$, respectively.
The complete split graph $K_{k,n-k}^*$  arises
from the complete bipartite graph $K_{k,n-k}$ by adding edges to make the
vertices in the partite set of size $k$ pairwise adjacent. Let $ {\mathcal G}_{n,p,k}$ be the family of complements of graphs  consisting of an $(n-k-1)$-regular
graph on $p$ vertices together with $n-p$ isolated vertices. We also can describe  $ {\mathcal G}_{n,p,k}$  as the family of $n$-vertex
graphs obtained from $K_n$ by deleting the edges of an $(n-k-1)$-regular graph on $p$ vertices.

In this paper we further investigate properties of the $R'$ index with regard to minimum degree $k$. We  found
in $G(k,n)$ graphs for which the variation of the Randi\'c index
 attains its minimum value.
 When $k \leq \frac n2$  the extremal graphs are complete split graphs $K_{k,n-k}^*$.
% , which   have vertices of only two degrees, i.e. degree $k$ and degree $n-1$, and the number of vertices of degree $k$  is $n-k$, while  the number of vertices of degree $n-1$ is $k$.
For $k \geq \frac n2$  the extremal graphs belong to the family ${\mathcal G}_{n,\frac n2,k}$.
% have vertices of two degrees $k$ and $n-1$, and the number of vertices of degree $k$ is $\frac n2$.
We proved next Theorem which match conjecture given by Aouchiche and Hansen about the Randi\'c index in [2].
\vskip\baselineskip
\noindent{\bf Theorem 1.} If $G$ is a graph of order $n$
 with $\delta(G) \geq k$, then
 $$
 R'(G) \geq \left \{
 \begin{array}{ll}
\frac n2 -\frac 12(\frac 1{ k} - \frac 1{ {n-1}}) k(n-k)     &\mbox  {if}\  k\leq \frac n2,   \\[0.4ex]
\frac n2 -\frac 12(\frac 1{ k} - \frac 1{ {n-1}}) \frac {n^2}{4}     &\mbox  {if}\   \frac n2  \leq k \leq n-2,
\end{array}
\right.
$$
For $ k \leq \frac n2$ equality holds if and only if $G =K_{k,n-k}^*$. For $ k \geq \frac n2$ equality holds
 if $n\equiv0 (\rm mod~ 4)$, or if $n\equiv 2 (\rm mod~ 4)$  and  $k$ is
odd, and
$G\in {\mathcal G}_{n,n/2,k}$.
\vskip\baselineskip
The proof is based on the approach  first time introduced in  [12].
%  and further reused in [11] and [14].

\vskip\baselineskip
\noindent {\bf 2. A QUADRATIC PROGRAMMING MODEL OF THE PROBLEM }
\vskip\baselineskip

 First, we will give some linear equalities and nonlinear
inequalities which must be satisfied in any
graph from the class $G(k,n)$. Let $x_{i,j}$ denote  the number of edges joining
vertices of degrees $i$ and $j$ and $n_i$ denote the number of vertices of degree $n_i$.
 The mathematical  description of
the problem $P$  to determine minimum of ~ $R'(G) =\sum_{k\leq i\leq j\leq  n-1} \frac {x_{i,j}} {\max \{ij\}} =\sum_{k\leq i \leq j\leq n-1} \frac {x_{i,j}} {j} $ ~ is:

$$  \min  \sum_{k\leq i \leq j  \leq n-1} \frac {x_{i,j}} {j}
$$
 subject to:
 $$
\begin{array}{c@{\:+\:}c@{\:+\:}c@{\:+\:}c@{\:+\:}c@{\;=\;}c}
2x_{k,k}  &  x_{k,k+1}   &   x_{k,k+2}   & \cdots      & x_{k,n-1}   &  kn_k,          \\
x_{k,k+1} & 2x_{k+1,k+1} &  x_{k+1,k+2}  & \cdots      & x_{k+1,n-1} &  (k+1)n_{k+1},  \\
% x_{k,k+2} &  x_{k+1,k+2} &  2x_{k+2,k+2} & \cdots      & x_{k+2,n-1} &   (k+2)n_{k+2},  \\
\multicolumn{6}{c}{\dotfill}                                                           \\
x_{k,n-1} &  x_{k+1,n-1} &   x_{k+2,n-1} & \cdots      & 2x_{n-1,n-1}&  (n-1)n_{n-1},    \\
\end{array}
\eqno{(1)}
$$
% and
$$
 n_k + n_{k+1} + n_{k+2} + \cdots + n_{n-1} = n,
 \eqno{(2)}
$$
$$
x_{i,j} \leq n_i n_j, \qquad \mbox{for}\quad k\leq i \leq n-1,
\quad i<j \leq n-1,     \eqno{(3)}
$$
% and
$$
x_{i,i} \leq {n_i \choose 2}, \qquad \mbox{for} \quad k \leq i \leq n-1,       \eqno{(4)}
$$
$$
x_{i,j},\  n_i \quad \mbox{ are non-negative integers, for }\ k\leq i \leq j \leq n-1.  \eqno{(5)}
$$
$(1-5)$ define a nonlinearly constrained optimization problem.

As it was done in [5], we  divide the first equality from (1)  by $k$, second by $k+1$, third by $k+2$ and so on, the last by $n-1$ and sum  them all, and  get
$$
\sum_{k\leq i\leq j\leq n-1}\left(\frac 1i + \frac 1j\right) x_{i,j} = n_k +n_{k+1} + n_{k+2} + \cdots + n_{n-1} = n,
$$
because of $(2)$. On the other side, $\frac 1j = \frac 12(\frac 1i + \frac 1j) - \frac 12(\frac 1i - \frac 1j)$. Then
\begin{eqnarray*}
R'(G) &  = &\sum_{k\leq i \leq j \leq n-1} \frac {x_{i,j}}{ {j}} = \frac 12 \sum_{k \leq i \leq j \leq n-1} \left( \frac 1i + \frac 1j - \left(\frac 1{ i} - \frac 1{ j}\right)\right) x_{i,j} \\
& = & \frac 12 \sum_{k \leq i \leq j \leq n-1} \left( \frac 1i + \frac 1j\right)x_{i,j} - \frac 12 \sum_{k \leq i \leq j \leq n-1} \left( \frac 1{i} - \frac 1{j}\right) x_{i,j}  \\
& = &  \frac n2 - \frac 12 \sum_{k \leq i \leq j \leq n-1} \left( \frac 1{i} - \frac 1{j}\right) x_{i,j}.
\end{eqnarray*}

We will  henceforth use the next expression (6) for the variation of the Randi\'c index:
$$
R'(G) = \frac n2 - \frac 12
  \sum_{k\leq i \leq j\leq n-1}
%  \sum_{k\leq i \leq n-1 \atop i \leq j \leq n-1 }
\left(\frac 1{i} - \frac 1{j}\right) x_{i,j}. \eqno {(6)}
$$

\par Define the function
$$
 \gamma = \sum_{k\leq i \leq j\leq n-1}
\left(\frac 1{ i} - \frac 1{ j}\right)
x_{i,j}.  \eqno {(7)}
$$
 Henceforth we will
consider the problem of maximizing  $\gamma $
instead of minimizing $ R'(G) $.
\vskip\baselineskip
\noindent  {\bf 3. PROOF OF THE FIRST PART OF THEOREM 1  $(k\leq \frac n2)$}
\vskip\baselineskip
\par\noindent{\bf PROOF:}  Since the minimum degree is $k$, it is evident that $n_{n-1} \leq k$.  Let $m$ be the index such that
$n_m + n_{m+1} + ...+ n_{n-2} + n_{n-1} \geq k$ and $n_{m+1} + ...+ n_{n-2} + n_{n-1} < k$.
We distinguish two subcases: (a) for such $m$ \ $n_m + ...+ n_{n-2} + n_{n-1} =k$, and (b) $n_m  + ...+ n_{n-2} + n_{n-1} >  k$.
\vskip\baselineskip
\par\noindent{\bf Subcase a.}~ $n_m + ...+ n_{n-2} + n_{n-1} =k$. We have:
\begin{eqnarray*}
\gamma   & =    & \sum_{k\leq i < j\leq n-1}\left(\frac{1}{{i}}-\frac{1}{{j}}\right) x_{i,j} = \sum_{j=k+1}^{n-1}\left(\frac{1}{{k}}-\frac{1}{{j}}\right) x_{k,j} \\
         & +    &  \sum_{j=k+2}^{n-1}\left(\frac{1}{{k+1}}-\frac{1}{{j}}\right) x_{k+1,j} +  \sum_{j=k+3}^{n-1}\left(\frac{1}{{k+2}}-\frac{1}{{j}}\right) x_{k+2,j}\\
         & +    &    \cdots +  \sum_{j=m}^{n-1}\left(\frac{1}{{m-1}}-\frac{1}{{j}}\right) x_{m-1,j} + \sum_{m\leq i < j\leq n-1}\left(\frac{1}{{i}}-\frac{1}{{j}}\right) x_{i,j}.
\end{eqnarray*}
$\sum_{j=i+1}^{n-1}(\frac{1}{{i}}-\frac{1}{{j}}) x_{i,j}$ represents weights of all edges which join  vertices of degree $i$, with vertices of degree $j$, \ $i+1\leq j\leq n-1$.
We give the maximum possible weights to these edges. Since $n_m + n_{m+1} +  \cdots + n_{n-1} = k$ and $\sum_{j=i+1}^{n-1}x_{i,j} \leq in_i$, first we join a vertex of degree $i$ to all $k$
vertices of degrees $n-1,\dots,m$ (maximum weights) and with $i-k$ vertices of other degrees $j$, \ $i+1 \leq j \leq m-1$. We will maximize the weights of these  last $i-k$ edges joining
a vertex of degree $i$ to $i-k$ vertices of degree $m-1$. Thus,
$$
\sum_{j=i+1}^{n-1}(\frac{1}{{i}}-\frac{1}{{j}}) x_{i,j} \leq n_i\left( \sum_{j=m}^{n-1}\left(\frac{1}{{i}}-\frac{1}{{j}}\right)n_j  + \left(\frac{1}{{i}}-\frac{1}{{m-1}}\right)(i-k)\right).
$$
Then
\begin{eqnarray*}
\gamma   & \leq & \sum_{k\leq i\leq m-1}g(i)n_i + \sum_{m\leq i < j\leq n-1}\left(\frac{1}{{i}}-\frac{1}{{j}}\right)x_{i,j} \\
         & =    &\Sigma_1 + \Sigma_2 .
\end{eqnarray*}
where $g(i) = \sum_{j=m}^{n-1}(\frac{1}{{i}}-\frac{1}{{j}})n_j  + (\frac{1}{{i}}-\frac{1}{{m-1}})(i-k)$.
Since $f(x) = x(\frac{1}{{x}}-\frac{1}{{y}})$, for $0<x< y$, is  a decreasing function, we have for $k+1 \leq i \leq m-1$,\ $m \leq j \leq n-1$:
$$
i(\frac{1}{{i}}-\frac{1}{{j}}) \leq k(\frac{1}{{k}}-\frac{1}{{j}}).
$$
Therefore
\begin{eqnarray*}
g(i) & \leq  & \frac ki\left(\sum_{j=m}^{n-1}(\frac{1}{{k}}-\frac{1}{{j}})n_j  + (\frac{1}{{k}}-\frac{1}{{m-1}})(i-k)\right) \\
     & =     & (1-\frac {i-k}i)\left(\sum_{j=m}^{n-1}(\frac{1}{{k}}-\frac{1}{{j}})n_j\right)  + \frac {k( i-k)}{i}(\frac{1}{{k}}-\frac{1}{{m-1}}) \\
     & =     & \sum_{j=m}^{n-1}(\frac{1}{{k}}-\frac{1}{{j}})n_j + \frac {i-k}i\left( (\frac{1}{{k}}-\frac{1}{{m-1}}) k - \sum_{j=m}^{n-1}(\frac{1}{{k}}-\frac{1}{{j}})n_j\right)\\
     & \leq  & \sum_{j=m}^{n-1}(\frac{1}{{k}}-\frac{1}{{j}})n_j,
\end{eqnarray*}
because $\sum_{j=m}^{n-1}n_j=k$ and $m \leq j \leq n-1$. Since $n_k + \dots + n_{m-1} = n-k$, we have
\begin{eqnarray*}
\Sigma_1   & = & \sum_{k\leq i\leq m-1}g(i)n_i \leq \left( \sum_{j=m}^{n-1}(\frac{1}{{k}}-\frac{1}{{j}})n_j\right) \sum_{i=k}^{m-1}n_i = (n-k) \sum_{j=m}^{n-1}(\frac{1}{{k}}-\frac{1}{{j}})n_j \\
           & = & \sum_{j=m}^{n-1}\left(\frac{1}{{k}}-\frac{1}{{n-1}}\right)(n-k)n_j  - \sum_{j=m}^{n-2}\left(\frac{1}{{k}}-\frac{1}{{n-1}}\right)(n-k)n_j  \\
           & + &  \sum_{j=m}^{n-2}\left(\frac{1}{{k}}-\frac{1}{{j}}\right)(n-k)n_j = \left(\frac{1}{{k}}-\frac{1}{{n-1}}\right)(n-k)k  \\
           & + & \sum_{j=m}^{n-2}\left(\left(\frac{1}{{k}}-\frac{1}{{j}}\right) -\left(\frac{1}{{k}}-\frac{1}{{n-1}}\right)\right)(n-k)n_j \\
           & = & \left(\frac{1}{{k}}-\frac{1}{{n-1}}\right)(n-k)k  - \sum_{j=m}^{n-2}\left(\frac{1}{{j}}-\frac{1}{{n-1}}\right)(n-k)n_j. \\
\end{eqnarray*}
Since $x_{i,j} \leq n_in_j$,\ $m\leq i < j\leq n-1$, and  $n_{n-1}= k- \sum_{j=m}^{n-2}n_j$, we have:
\begin{eqnarray*}
\Sigma_2 & = & \sum_{m\leq i < j\leq n-1}\left(\frac{1}{{i}}-\frac{1}{{j}}\right)x_{i,j} \leq \sum_{m\leq i < j\leq n-1}\left(\frac{1}{{i}}-\frac{1}{{j}}\right)n_in_j \\
         & = &  \sum_{i=m}^{n-2}\left(\frac{1}{{i}}-\frac{1}{{n-1}}\right)n_i \left (k - \sum_{j=m}^{n-2}n_j \right) + \sum_{m\leq i < j\leq n-2}\left(\frac{1}{{i}}-\frac{1}{{j}}\right)n_in_j \\
         & = & k \sum_{i=m}^{n-2}\left(\frac{1}{{i}}-\frac{1}{{n-1}}\right)n_i - \sum_{i=m}^{n-2}\left(\frac{1}{{i}}-\frac{1}{{n-1}}\right)n_i^2 \\
         & + &  \sum_{m\leq i < j\leq n-2}\left( (\frac{1}{{i}}-\frac{1}{{j}}) - (\frac{1}{{i}}-\frac{1}{{n-1}})  - (\frac{1}{{j}}-\frac{1}{{n-1}}) \right) n_in_j \\
         & = & k \sum_{i=m}^{n-2}\left(\frac{1}{{i}}-\frac{1}{{n-1}}\right)n_i -  \sum_{i=m}^{n-2}\left(\frac{1}{{i}}-\frac{1}{{n-1}}\right)n_i^2 \\
         & - &  2 \sum_{m\leq i < j\leq n-2}\left(\frac{1}{{j}}-\frac{1}{{n-1}}\right)n_in_j\\
%         & = &  k\sum_{i=m}^{n-2}\left(\frac{1}{{i}}-\frac{1}{{n-1}}\right)n_i - \sum_{i=m}^{n-2}\left(\frac{1}{{i}}-\frac{1}{{n-1}}\right)n_i \sum_{j=m}^{n-2}n_j  \\
%        & + &  \sum_{m\leq i < j\leq n-2}\left(\frac{1}{{i}}-\frac{1}{{j}}\right)n_in_j  =  k\sum_{i=m}^{n-2}\left(\frac{1}{{i}}-\frac{1}{{n-1}}\right)n_i\\
%       & + &  \sum_{i=m}^{n-2}n_i \left( (\frac 1i - \frac 1{i+1})n_{i+1} + (\frac 1i - \frac 1{i+2})n_{i+2} + \dots + (\frac 1i - \frac 1{n-2})n_{n-2}\right. \\
%         & - &  \left.(\frac 1i - \frac 1{n-1})(n_m +n_{m+1} + \dots + n_i + \dots + n_{n-2})\right) \\
%        & = & k\sum_{i=m}^{n-2}\left(\frac{1}{{i}}-\frac{1}{{n-1}}\right)n_i - \sum_{i=m}^{n-2}n_i \left(\frac 1i - \frac 1{n-1}\right)\left(\sum_{j=m}^{i}n_j\right)\\
%        & + & \sum_{i=m}^{n-2}n_i\sum_{j=i+1}^{n-2}\left( (\frac 1i - \frac 1j) - (\frac 1i - \frac 1{n-1})\right)n_j \\
 %        & = & k\sum_{i=m}^{n-2}\left(\frac{1}{{i}}-\frac{1}{{n-1}}\right)n_i - \sum_{i=m}^{n-2}n_i \left(\frac 1i - \frac 1{n-1}\right)\left(\sum_{j=m}^{i}n_j\right)\\
%         & - & \sum_{i=m}^{n-2}n_i\sum_{j=i+1}^{n-2}\left( \frac 1j - \frac 1{n-1}\right)n_j. \\
\end{eqnarray*}
Thus,
\begin{eqnarray*}
\gamma & \leq & \Sigma_1  +  \Sigma_2 \leq \left(\frac{1}{{k}}-\frac{1}{{n-1}}\right)(n-k)k    \\
       & -    & \sum_{i=m}^{n-2}\left(\frac{1}{{i}}-\frac{1}{{n-1}}\right)(n-k)n_i \\
       & +    & k\sum_{i=m}^{n-2}\left(\frac{1}{{i}}-\frac{1}{{n-1}}\right)n_i -   \sum_{i=m}^{n-2}\left(\frac{1}{{i}}-\frac{1}{{n-1}}\right)n_i^2 \\
       & -    &  2 \sum_{m\leq i < j\leq n-2}\left(\frac{1}{{j}}-\frac{1}{{n-1}}\right)n_in_j  = \left(\frac{1}{{k}}-\frac{1}{{n-1}}\right)(n-k)k \\
       &  -   &  \sum_{i=m}^{n-2}\left(\frac{1}{{i}}-\frac{1}{{n-1}}\right)(n-2k)n_i  -  \sum_{i=m}^{n-2}\left(\frac{1}{{i}}-\frac{1}{{n-1}}\right)n_i^2    \\
       & -    &  2 \sum_{m\leq i < j\leq n-2}\left(\frac{1}{{j}}-\frac{1}{{n-1}}\right)n_in_j  \leq \left(\frac{1}{{k}}-\frac{1}{{n-1}}\right)(n-k)k.    \\
%        \sum_{i=m}^{n-2}n_i \left(\frac 1i - \frac 1{n-1}\right)\left(\sum_{j=m}^{i}n_j\right)\\
 %      & -    &  \sum_{i=m}^{n-2}n_i\sum_{j=i+1}^{n-2}\left( \frac 1j - \frac 1{n-1}\right)n_j \leq \left(\frac{1}{{k}}-\frac{1}{{n-1}}\right)(n-k)k. \\
\end{eqnarray*}
The last inequality follows because $k \leq \frac n2$.
     Equality holds when $n_i=0$ for $k+1 \leq i \leq n-2$, $n_k = n-k$,\ $n_{n-1} = k$,\  $x_{k,n-1} = (n-k)k $,\  $x_{n-1,n-1} =  {k \choose 2}$,
     and all other $x_{i,j}$ are equal to zero.
     Thus, graphs for which variation of the Randi\'c index attains its minimum value are $K_{k,n-k}^*$.
\vskip\baselineskip
\par\noindent{\bf Subcase b.}~  We put $n_m = n_{m'} + n_{m''}$, such that $n_{m''} + n_{m+1} + ... + n_{n-1} =k$. Then $n_k + \cdots + n_{m-1} + n_{m'} = n-k$. We will color the vertices of degree $m$ with red and white, such that the number of red vertices is $n_{m''}$. Denote by $x_{i,m'} (x_{i,m''})$  for $i\neq m$, the number of edges between vertices of degree $i$ and the white (red) vertices
of degree $m$, by  $x_{m',m'}$ $(x_{m'',m''})$ the number of edges between white (red) vertices of degree $m$, and by $x_{m',m''}$ the number of edges between  white and red vertices of degree $m$. Then $x_{i,m} =
x_{i,m'} + x_{i,m''}$  for $i\neq m$, and $x_{m,m} = x_{m',m'} + x_{m',m''} + x_{m'',m''}$. We will replace system $(1)$ by:
$$
\begin{array}{c@{\:+\:}c@{\:+\:}c@{\:+\:}c@{\:+\:}c@{\:+\:}c@{\:+\:}c@{\:+\:}c@{\:=\:}l }
x_{k,i}   &  \cdots &  x_{i,m-1}    &  x_{i,m'}   &  x_{i,m''}    &  x_{i,m+1}   &  \cdots &   x_{i,n-1}   &  in_i, \\
\multicolumn{9}{r}{ k\leq i\leq n-1, i\neq m,}                                                                      \\
x_{k,m'}  &  \cdots &  x_{m-1,m'}   &  2x_{m',m'} &  x_{m',m''}   &  x_{m',m+1}  &  \cdots &   x_{m',n-1}  &  mn_{m'}, \\
x_{k,m''} &  \cdots &  x_{m-1,m''}  &  x_{m',m''} &  2x_{m'',m''} &  x_{m'',m+1} &  \cdots &   x_{m'',n-1} &  mn_{m''}, \\
\end{array} \eqno (\tilde 1)
$$
We will proceed similarly  as in the subcase a. The rest of the proof is omitted, because it is similar to the one of subcase a.
\hfill  $\Box$
\vskip\baselineskip
\vskip\baselineskip
\noindent   {\bf 4. PROOF OF THE SECOND PART OF THEOREM 1 $(k\geq \frac n2)$}
\vskip\baselineskip

 We put:
$$
\begin{array}{lcl}
  x_{i,j} = n_i n_j - y_{i,j} \qquad           & \mbox {for}    & \quad k \leq i \leq n-1, \ i < j \leq n-1,  \\
  x_{i,i} =  {n_i \choose 2} - y_{i,i} \qquad  & \mbox  {for}   &\quad k \leq i \leq n-1.
\end{array}  \eqno{(8)}
$$
A vertex of degree $n-1$ is adjacent to all other vertices. Thus $y_{i,n-1} = 0$\ for $k\leq i \leq n-1$ and $n_{n-1} \leq k$,
or the minimum degree would be greater than $k$.
After substitution of $x_{i,j}$ and $x_{i,i}$
from $(8)$ into the  function  $\gamma $  and $(1)$, we rewrite the optimization problem using the same objective function (call the rewritten problem $\overline P$) as:
$$
\max \sum_{k\leq i < j\leq n-1}
\left(\frac 1{i} - \frac 1{j}\right)n_i n_j -
\sum_{k\leq i < j\leq n-2}
\left(\frac 1{i} - \frac 1{j}\right)y_{i,j}
$$
subject to
$$
\begin{array}{c@{\:+\:}c@{\:+\:}c@{\:+\:}c@{\:+\:}c@{\;=\;}l}
2y_{k,k}  &  y_{k,k+1}   &   y_{k,k+2}   & \cdots      & y_{k,n-2}   &  (n-k-1)n_k,      \\
y_{k,k+1} & 2y_{k+1,k+1} &  y_{k+1,k+2}  & \cdots      & y_{k+1,n-2} &  (n-k-2)n_{k+1},  \\
% y_{k,k+2} &  y_{k+1,k+2} &  2y_{k+2,k+2} & \cdots      & y_{k+2,n-2} &  (n-k-3)n_{k+2},  \\
\multicolumn{6}{c}{\dotfill}                                                           \\
y_{k,n-2} &  y_{k+1,n-2} &   y_{k+2,n-2} & \cdots      & 2y_{n-2,n-2}&    n_{n-2},    \\
\end{array}
\eqno{(1')}
$$
$$
 n_k + n_{k+1} + n_{k+2} + \cdots + n_{n-1} = n,
 \eqno{(2)}
$$
$$
\begin{array}{lclr}
\phantom{0000000} n_ i\geq 0,       & \phantom{000} \mbox {for}   &       k \leq i \leq n-1,                    & \phantom{000000000000000} (9)      \\
\phantom{0000000} y_{i,j} \geq 0,   & \phantom{000}  \mbox {for}  &       k \leq i \leq n-2,\ i \leq j \leq  n-2,   &  (10)    \\
% \phantom{0000000} y_{i,i} \geq 0,   &  \phantom{000}  \mbox {for} &       k \leq i \leq n-2,                    &  (10)    \\
\phantom{0000000} n_{n-1} \leq k,   &                &                                             &  (11)
\end{array}
$$
$$
y_{i,j},\  n_i \quad \mbox{ are  integers  for }\ k\leq i \leq j \leq n-1.  \eqno{(5')}
$$
We obtained equalities $(1')$ from the corresponding equalities $(1)$.
Let $(n_k,n_{k+1},\\ \ldots,n_{n-1},y_{k,k},y_{k,k+1},\ldots,y_{n-2,n-2})$  be a  feasible point for  $\overline P$; we use $\Omega$ or $ ({N,Y})$
to denote this point. Let  $ \gamma_1 = \sum_{k\leq i < j  \leq n-1}(\frac 1{i} - \frac 1{
j})n_i n_j $  and  $ \gamma_2 = -\sum_{k\leq i < j \leq n-2} (\frac 1{ i} - \frac 1{j})
y_{i,j} $. Now $ \max \gamma \leq \max\gamma_1 + \max\gamma_2$, where the maxima are subject to
$(1'), (2), (9-11), (5') $. It is evident that
$\max \gamma_2 = 0$, and it is achieved by setting $y_{i,j} = 0$ for
$k \leq i \leq n-2$   and  $i < j \leq n-2$ and  setting $y_{i,i} = \frac {(n-i-1)n_i}2$  for  $ k \leq i \leq n-2$.
The variables $n_i$ must satisfy $(2), (9)$, $(11)$ and  $(5')$.
Hence, there are many extreme points for $\gamma _2$.
Let us denote by $(n_k^*,n_{k+1}^*,\ldots,n_{n-1}^*)$  or $N^*$ the optimal point for
$\gamma_1$. Let $Y^* = (y_{k,k}^*,y_{k,k+1}^*,\ldots,y_{n-2,n-2}^*)$,
where \ $y_{i,j}^* = 0$ for  $ i \neq j $ and $y_{i,i}^* = \frac {(n-i-1)n_i^*}2$. Note that $Y^*$  is the optimal point for $\gamma_2$
if $y^*_{i,j}$ are integers,
and  $(N^*,Y^*)$
% which we denote $\Omega^*$,
  will be the optimal point for $\gamma$.
  In order to find $N^*$ we can neglect
 constraints $(1')$  and $(10)$, because for $\gamma_1$ only constraints $(2), (9)$ and $(11)$
are relevant. We omit constraint (11), because it is not
necessary and would complicate the calculation. We also neglect constraint $(5')$, but we will keep it in mind.
% and introduce new constraints when it is necessary.

We will need the following theorems.

\vskip\baselineskip
\noindent {\bf Theorem 1.4.10 from [ 15].} A two times differentiable function $f$ on open convex set $C$ is concave if and only if Hessian matrix
$$
 H(x) = \left[ \frac {\partial^2 f(x)}{\partial x_i \partial x_j}\right]
 $$
 is negative-semidefinite matrix for $\forall x\in C$.
 \vskip\baselineskip
\noindent {\bf Generalized Sylvester's cryterion.} A $n\times n$ Hermitian matrix $A =(a_{i,j})$   is negative-definite if and only if
members of the sequence $1,D_1,D_2,\dots,D_n$ change the sign, where $D_i$ are the principal minors, that is ($D_1<0,D_2>0,\dots)$.
\vskip\baselineskip
 From (2), we have $n_{n-1} = n- \sum_{j=k}^{n-2}n_j$.
We rewrite $\gamma_1$:
\begin{eqnarray*}
\gamma_1 & = & \sum_{k\leq i < j\leq n-2}\left(\frac{1}{{i}}-\frac{1}{{j}}\right)n_in_j
 + \sum_{i=k}^{n-2}\left(\frac{1}{{i}}-\frac{1}{{n-1}}\right)n_i(n-\sum_{j=k}^{n-2} n_j) \\
& = & n \sum_{i=k}^{n-2}\left(\frac{1}{{i}}-\frac{1}{{n-1}}\right)n_i -
\sum_{i=k}^{n-2}\left(\frac{1}{{i}}-\frac{1}{{n-1}}\right)n_i^2 \\
& + &
\sum_{k\leq i < j\leq n-2}\left( (\frac{1}{{i}}-\frac{1}{{j}}) - (\frac{1}{{i}}-\frac{1}{{n-1}})  -
(\frac{1}{{j}}-\frac{1}{{n-1}}) \right) n_in_j \\
&=& n \sum_{i=k}^{n-2}\left(\frac{1}{{i}}-\frac{1}{{n-1}}\right)n_i -
\sum_{i=k}^{n-2}\left(\frac{1}{{i}}-\frac{1}{{n-1}}\right)n_i^2 \\
& - &  2 \sum_{k\leq i < j\leq n-2}\left(\frac{1}{{j}}-\frac{1}{{n-1}}\right)n_in_j  ~~~~~~~~~~~~~~~~~~~~~~~~~~~~~~~~~~~~~~~~(12)\\
\end{eqnarray*}
Define a function  $\bar \gamma_1$ by $ \bar \gamma_1(n_k,...,n_{n-2}) = n \sum_{i=k}^{n-2}\left(\frac{1}{{i}}-\frac{1}{{n-1}}\right)n_i -
\sum_{i=k}^{n-2}\left(\frac{1}{{i}}-\frac{1}{{n-1}}\right)n_i^2  -
2 \sum_{k\leq i < j\leq n-2}\left(\frac{1}{{j}}-\frac{1}{{n-1}}\right)n_in_j$
(see (12)).  Let $X = \{ (n_k,...,n_{n-1}) \mid n_k + ... + n_{n-1} =n\}$.
 Note that $ \gamma_1(n_k,...,n_{n-1}) =  \bar \gamma_1(n_k,...,n_{n-2})$ for $(n_k,...,n_{n-1})\in X$.
  We will study $\bar \gamma_1$  on  ${\mathbb R}^{n-k-1}$ instead of $\gamma_1$ on $X$.
  The  point $(n_{k},...,n_{n-2})\in {\mathbb R}^{n-k-1}$
  corresponds to  $( n_{k},...,n_{n-2},n-\sum_{j=k}^{n-2}n_j)\in {\mathbb R}^{n-k}$ on the set $X$.
Let us notice that
$\bar \gamma_1$ on ${\mathbb R}^{n-k-1}$ is concave function. The  $j$-th  principal minor  is
$$
D_j = (-2)^{j}\left |
\begin{array}{ccccc}
(\frac 1k - \frac 1{n-1})       & (\frac 1{k+1} - \frac 1{n-1})    & (\frac 1{k+2} - \frac 1{n-1})   &   \dots     & (\frac 1{k+j-1} - \frac 1{n-1})   \\
(\frac 1{k+1} - \frac 1{n-1})   & (\frac 1{k+1} - \frac 1{n-1})    & (\frac 1{k+2} - \frac 1{n-1})   &  \dots      &  (\frac 1{k+j-1} - \frac 1{n-1})   \\
(\frac 1{k+2} - \frac 1{n-1})   & (\frac 1{k+2} - \frac 1{n-1})    & (\frac 1{k+2} - \frac 1{n-1})   &  \dots      &  (\frac 1{k+j-1} - \frac 1{n-1})   \\
\vdots                          &                       \vdots     & \vdots                          & \ddots      & \vdots                              \\
(\frac 1{k+j-1} - \frac 1{n-1}) & (\frac 1{k+j-1} - \frac 1{n-1})  & (\frac 1{k+j-1} - \frac 1{n-1}) & \dots       & (\frac 1{k+j-1} - \frac 1{n-1})   \\
\end{array} \right|.
$$
It is not difficult to find that $D_j = 2^j(-1)^j(\frac 1k - \frac 1{k+1})(\frac 1{k+1} - \frac 1{k+2})(\frac 1{k+2} - \frac 1{k+3})\cdots (\frac 1{k+j-2} - \frac 1{k+j-1})(\frac 1{k+j-1} - \frac 1{n-1})$.
Using Sylvester's cryterion we conclude that $\bar \gamma_1$ is concave function.

%%%%%%%%%%%%%%%%%%%%%%%%%%%%%%%%%%%%%%%%%%%%%%%%%%%%%%%%%%%%%%%%%%%%%%%%%%%%%%%%%%%%%%%%%%%%%%%%%%%%%
%\vskip\baselineskip
%\noindent {\bf  5.  CASE 1: ${\bf n\equiv 0( \rm \bf mod\ 4)}$, OR ${\bf n\equiv 2(\rm \bf mod\ 4)}$ AND  $\bf k$ IS ODD}
%\vskip\baselineskip
\par

 We consider the  problem $\overline{P^1}$ of maximizing $\bar \gamma_1$:
$$
\max ~ n \sum_{i=k}^{n-2}\left(\frac{1}{{i}}-\frac{1}{{n-1}}\right)n_i  - \sum_{i=k}^{n-2}\left(\frac{1}{{i}}-\frac{1}{{n-1}}\right)n_i^2
 -   2 \sum_{k\leq i < j\leq n-2}\left(\frac{1}{{j}}-\frac{1}{{n-1}}\right)n_in_j
%\sum_{k\leq i < j\leq n-1}\left(\frac{1}{{i}}-\frac{1}{{j}}\right)n_in_j,
$$
subject to
% $$
% n_k+n_{k+1}+\cdots+n_{n-1}=n,     \eqno{(2)}
% $$
$$
n_i \geq 0~~\mbox{for}~~k\leq i \leq n-2,      \eqno{(9)}
$$
instead of the problem: $\max \gamma_1$ subject to $(2)$ and $(9)$.
 Let $N$ denote a feasible point
$(n_k,\ldots ,n_{n-2})$
for problem $\overline{P^1}$.
We will show that   $N_1^*$ is an optimal point for the problem
$\overline {P^1}$, where $N_1^*$ is defined by
 $n_k=\frac n2, n_i=0$ for $k+1\leq i\leq n-2$.
\vskip\baselineskip
\noindent {\bf Lemma 1.} {\it The function $\gamma_1$, subject to (2) and (9), attains its
maximum   value $\gamma_1^*$ equal to
$\frac{n^2}{4}(\frac{1}{{k}}-\frac{1}{{n-1}})$
at the point $(\frac n2,0,0,\ldots,0,\frac n2)\in {\mathbb R}^{n-k}$.}
\vskip\baselineskip
\par\noindent{\bf PROOF:}
We distinguish two subcases: (1a) $\Delta(G)=n-1$, and (1b) $\Delta(G)<n-1$.
\vskip\baselineskip
\noindent {\bf Subcase 1a.}~ We will find point $N=(n_{k},...,n_{n-2})$ for which $\partial {\bar \gamma_1} / \partial n_i =0$,  $k\leq i\leq n-2$,
respectively:
$$
n \left(\frac{1}{{k}}-\frac{1}{{n-1}}\right)  - 2\left(\frac{1}{{k}}-\frac{1}{{n-1}}\right)n_k
 -   2 \sum_{k < j\leq n-2}\left(\frac{1}{{j}}-\frac{1}{{n-1}}\right)n_j    = 0,  \eqno(13)
 $$
\begin{eqnarray*}
&& n \left(\frac{1}{{i}}-\frac{1}{{n-1}}\right)  - 2\left(\frac{1}{{i}}-\frac{1}{{n-1}}\right)n_i
 -   2 \sum_{k \leq j < i}\left(\frac{1}{{i}}-\frac{1}{{n-1}}\right)n_j    \\
&  - &   2 \sum_{i < j \leq n-2}\left(\frac{1}{{j}}-\frac{1}{{n-1}}\right)n_j = 0, \qquad\quad
 \mbox{for}~~k+1\leq i\leq n-2, \qquad\quad\   (14)
\end{eqnarray*}
% and
% $$
% \lambda_in_i=0~~~\mbox{for}~~~k\leq i\leq n-2.     \eqno{(15)}
% $$
It is easy to see that the point $(\frac n2,0,0,\ldots,0)$   satisfies equalities
conditions $(13-14)$. Since $\bar \gamma_1$ is a concave function the point $(\frac n2,0,0,\ldots,0)$ is a maximum point and the maximum value $\bar \gamma_1^* =\gamma_1^* $
is
$$
\frac{n^2}{4}(\frac{1}{{k}}-\frac{1}{{n-1}}).
$$
To point $N_1^* = (\frac n2,0,0,\ldots,0)\in {\mathbb R}^{n-k-1}$ corresponds point $ N^* =(\frac n2,0,0,\ldots,0,\frac n2)\in {\mathbb R}^{n-k}$.
\vskip\baselineskip
\noindent {\bf Subcase 1b.} Let $m =\Delta (G)$. For this case the proof is similar to that of subcase 1a
and  is omitted. The maximum value of
$\gamma_1$ is
$$
\gamma_1^m=\left(\frac{1}{{k}}-\frac{1}{{m}}\right)\frac{n^2}{4}.
$$
This value is attained at the point $n_k=\frac n2, n_i=0$ for
$k+1 \leq i \leq m-1$, and $n_m=\frac n2$.
Since $\gamma_1^* > \gamma_1^m$  for this case,  we conclude that
$\gamma_1^*$ is the maximum value, attained at  $N^*$ on the
set of all feasible points.  \hfill $\square$
\vskip\baselineskip
 We have proved that $\gamma_1$
attains its maximum at  $N^*$. Now observe  that $\gamma_2$ attains
its maximum value, which is 0, at the point $Y^*_1$ defined by setting $y_{k,k} = \frac {(n-k-1)n}4$
 and all other $y_{i,j} = 0$.
% These values are integers, and
We conclude that the variation of the
Randi\'c index attains its minimum on
graphs  for which $n_k=n_{n-1}=n/2, x_{k,k}=n(2k-n)/8,
x_{k,n-1}=n^2/4,x_{n-1,n-1}=n(n-2)/8$, and all other
$x_{i,j},x_{i,i}$ and $n_{i}$ are equal to zero. These values are  integers only if $n\equiv0 (\rm mod~ 4)$, or if $n\equiv 2 (\rm mod~ 4)$  and  $k$ is
odd, and  these graphs lie in ${\mathcal G}_{n,n/2,k}$.

  Thus, we come to our conclusion:
\vskip\baselineskip
\noindent {\bf The second part of Theorem 1.}~~{\it  If $G \in G(k,n)$, then
$$
R'(G)\geq
\frac{n}{2}-\frac{n^2}{8}\left(\frac{1}{{k}}-\frac{1}{{n-1}}\right).
$$
If
$n\equiv0 (\rm mod~ 4)$, or if $n\equiv 2 (\rm mod~ 4)$  and  $k$ is
odd, this value is attained by graphs
$G\in {\mathcal G}_{n,n/2,k}$ for which  $n_k=n_{n-1}=n/2,
x_{k,k}=n(2k-n)/8, x_{k,n-1}=n^2/4,x_{n-1,n-1}=n(n-2)/8$, and
all other $x_{i,j},x_{i,i}$ and $n_{i}$ are equal to
zero.
}
\vskip\baselineskip
Furthermore, if $G\in {\mathcal G}_{n,p,k}$, then  $R'(G) = \frac n2 -\frac 12(\frac 1{ k} - \frac 1{{n-1}}) p(n-p)$. Now we give conjecture about extremal graphs
for the other parity of $n$ and $k$.
\vskip\baselineskip
 \noindent{\bf Conjecture.} If $G \in G(k,n)$ and if
 $$
 p= \left \{
\begin{array}{ll}
% \frac n2                                & \mbox  {if}\ n\equiv 0 (\mbox {mod}\ 4);\quad n\equiv 2(\mbox {mod}\ 4),\   k\  \mbox  {is odd},\quad (\mbox {\bf it is done})\\[0.4ex]
 \lfloor \frac n2 \rfloor \ \mbox {or}\ \lceil \frac n2 \rceil
                                        & \mbox  {if}\ n\equiv 1 (\mbox {mod}\ 4),\ k\  \mbox  {is even};\
                                         \ n\equiv 3 (\mbox {mod}\ 4),\ k\  \mbox  {is even},\\[0.4ex]
 \lfloor \frac n2 \rfloor                & \mbox  {if}\ n\equiv 1 (\mbox {mod}\ 4),\ k\  \mbox  {is odd}, \\[0.4ex]
 \frac {n-2}2\  \mbox {or}\ \frac {n+2}2  & \mbox  {if}\ n\equiv 2(\mbox {mod}\ 4),\   k\  \mbox  {is even},   \\[0.4ex]
\lceil \frac n2 \rceil                  & \mbox  {if}\ n\equiv 3(\mbox {mod}\ 4),\  k\  \mbox  {is odd}, \\
 \end{array}\right.
 $$
 then
  $$
 R'(G) \geq   \frac n2 -\frac 12(\frac 1{k} - \frac 1{ {n-1}}) p(n-p) \quad\quad    \mbox  {if}\   \frac n2  < k \leq n-2,
 $$
 where $p$ and $n$ are given above. Equality holds if and only  if $G \in {\mathcal G}_{n,p,k}$.
 \vskip\baselineskip

 The proof of this conjecture is more complicated than for the case  $n\equiv0 (\rm mod~ 4)$, or if $n\equiv 2 (\rm mod~ 4)$  and  $k$ is
 odd  and we leave it as an open problem.
More information on $ \gamma_1$ function could be obtained using maximizing technique  given in [ 14] or approach in [11].
% and could be done in the  manner similar to that in [?].

Let  $G(k,m, n)$ be the  set of connected simple $n$-vertex graphs
with minimum  vertex degree $k$ and  maximum  vertex degree $m$, where $k \leq m\leq n-2$.
Let  ${\mathcal G}_{n,p,k,m}$ be the family of complements of graphs
%$G_{n,p,k,m}$
consisting of an $(n-k-1)$-regular
graph on $p$ vertices  and $(n-m-1)$-regular
graph on $n-p$ vertices. Since the proof of the next theorem is  similar to those  of Theorems 1, we  omit
them and just write  down the theorem.
\vskip\baselineskip
\noindent{\bf Theorem 2.} {\it If $G$ is a graph of order $n$
 with $\delta(G) \geq k$ and $\triangle(G) \leq m$, then
 $$
 R'(G) \geq \left \{
 \begin{array}{ll}
\frac n2 -\frac 12(\frac 1{ k} - \frac 1{ {m}}) k(n-k)     &\mbox  {if}\  k\leq \frac n2,\quad\mbox  {and}\  n-k \leq m\leq n-2    \\[0.4ex]
\frac n2 -\frac 12(\frac 1{ k} - \frac 1{ {m}}) \frac {n^2}{4}     &\mbox  {if}\   \frac n2  \leq k \leq m \leq n-2,
\end{array}
\right.
$$
 For $ k \leq \frac n2$, equality holds  if $k$ is even, or $k$ and $n-m$ are odd, and  this value is attained by graphs
for which  $n_k=n-k, n_{m}= k,
x_{k,m}=(n-k)k,x_{m,m}=k(k+m-n)/2$, and
all other $x_{i,j},x_{i,i}$ and $n_{i}$ are equal to
zero. For $ k \geq \frac n2$ equality holds if
$n\equiv0 (\rm mod~ 4)$, or if  $n\equiv 2 (\rm mod~ 4)$, and $k$ and $m$ are
odd,  and this value is attained by graphs $G\in {\mathcal G}_{n,n/2,k,m}$ for which   $n_k=n_{m}=n/2,
x_{k,k}=n(2k-n)/8, x_{k,m}=n^2/4,x_{m,m}=n(2m-n)/8$, and
all other $x_{i,j},x_{i,i}$ and $n_{i}$ are equal to
zero. }

\vskip\baselineskip

\vskip\baselineskip
\noindent {\bf Acknowledgement:} This research was supported by
Serbian Ministry for Education and Science, Project No. 174010 {\it "Mathematical Models and Optimization Methods for Large Scale Systems} and Project No.
174033 {\it "Graph Theory and Mathematical Programming with
Applications to Chemistry and Computing "}.
%\newpage

\vskip\baselineskip
\par\noindent {\bf References}
\vskip\baselineskip
\par\noindent [1] V. Andova, M. Knor, P. Poto\v cnik,  R. \v Skrekovski, {\it On a variation of the Randi\'c index},
Australian Journal of Combinatorics, 56 (2013), 61-67.

\par\noindent [2] M. Aouchiche, P. Hansen, {\it On a conjecture about the Randi\'c index},
Discrete Mathematics, 307 (2), 2007, 262-265.

\par\noindent [3 ] B. Bollob\'as, P. Erd\H os, {\it Graphs of Extremal Weights},
Ars Combinatoria,  50 (1998), 225 - 233.

\par\noindent [4 ] M. Cygan, M. Pilipczuk, R. \v Skrekovski,  {\it On the inequality between radius and Randi\'c index for Graphs}, MATCH - Commun. Math. Comput. Chem.,
67 (2012), 451-466.

\par\noindent [5 ] T. Divni\'c, Lj. Pavlovi\'c,  {\it Proof of the first part of the conjecture of Aouchiche and Hansen about the Randi\'c index},
{ Discrete Applied Mathematics}, Vol.161, Issues 7-8,(2013), pp. 953-960.

% \par\noindent [6 ]  T. Divni\'c,  Lj. Pavlovi\'c, B. Liu, {\it Extremal graphs for the Randi\'c  index when minimum, maximum degrees and order of graphs are odd},
% { Optimization}, Vol. 64, No. 9, 2015, pp. 2021-2038.

\par\noindent [6 ] Dvo\v{r}\'ak, B. Lidicky, R. \v Skrekovski, {\it Randi\'c index and the diameter of a graph},
European Journal of Combinatorics, 32 (2011),434-442.

\par\noindent [ 7 ] S. Fajtlowicz (1998), {\it Written on the Wall}, Conjectures derived
on the basis of the program Galatea Gabriella Graffiti, University of Houston.

 \par\noindent [ 8] L.B. Kier, L.H. Hall, {\it Molecular Connectivity in
 Chemistry and Drug Research}, Academic Press, New York (1976).

 \par\noindent [ 9 ] L.B. Kier, L.H. Hall, {\it Molecular Connectivity in
 Structure-Activity Analysis}, Research Studies Press-Wiley, Chichester (UK) (1986).

 \par\noindent [10 ] X. Li, I. Gutman, {\it Mathematical Aspects of Randi\'c-Type
  Molecular Structure Descriptors}, Mathematical Chemistry Monographs No. 1, Kragujevac,
 2006, pp. VI+330.

\par\noindent [11] B. Liu, Lj. Pavlovi\'c, J. Liu, T. Divni\'c, M. Stojanovi\'c, {\it On the conjecture of Aouchiche and Hansen about the Randi\'c index},
 Discrete Mathematics, Vol. 313, Issue 3 (2013), pp. 225-235.

\par\noindent [12]  Lj. Pavlovi\'c, T. Divni\'c, {\it A  quadratic programming approach
to the Randi\'c index}, European Journal of Operational
Research, 2007, Vol. 176, Issue 1, pp. 435-444.

\par\noindent [13] M. Randi\'c, {\it On characterization of molecular
branching}, J. Am. Chem. Soc., 97 (1975), 6609 -6615.

\par\noindent [14 ] I. Tomescu, R. Marinescu-Ghemeci, G. Mihai, {\it On dense graphs having minimum Randi\'c index},
Romanian J. of Information Science and Technology  (ROMJIST),
Vol. 12, (4), 2009,  pp. 455-465.

\par\noindent [15]  V. Vuj\v ci\'c, M. A\v si\'c, N. Mili\v ci\'c, {\it Matemati\v cko programiranje},
Matemati\v cki Institut, Beograd (1980).

\end{document}